# Apodictic discourse and the Cauchy-Bunyakovsky-Schwarz inequality

Satyanad Kichenassamy[1]

**Abstract.** Bunyakovsky's integral inequality (1859) is one of the familiar tools of modern Analysis. We try and understand what Bunyakovsky did, why he did it, why others did not follow the same path, and explore some of the mathematical (re)interpretations of his inequalities. This is achieved by treating the texts as discourses that provide motivation and proofs by their very discursive structure, in addition to what meets the eye at first reading. Bunyakovsky paper is an outgrowth of the mathematical theory of mean values in Cauchy's work (1821), but viewed from the point of view of Probability and Statistics. Liouville (1836) gave a result that implies Bunyakovsky's inequality, but did not identify it as significant because his interests lay elsewhere. Grassmann (1862) stated the inequality in abstract form but did not prove it for reasons that can be identified. Finally, by relating the result to quadratic binary forms, Schwarz (1885) opened the way to a geometric interpretation of the inequality that became important in the theory of integral equations. His argument is the source of one of the proofs most commonly taught nowadays. At about the same time, the Rogers-Hölder inequality suggested generalizations of Cauchy's and Bunyakovsky's results in an entirely different direction. Later extensions and reinterpretations show that no single result, even now, subsumes all known generalizations.

## 1. Introduction

Viktor Yakovlevich Bunyakovsky Викторъ Яковлевичъ Буняковскій (1804-1889) was one of the important figures of the Saint-Petersburg school of Mathematics in the nineteenth century, much respected and much loved[2]. He retired from the university in 1859[3] and, in the same year, published his most famous paper « Sur quelques inégalités concernant les intégrales ordinaires et les intégrales aux différences finies »[4], that contains the inequality, labeled (C) in his paper (p. 4),

$$\int_{x_0}^{X} \varphi(x)^2 dx \int_{x_0}^{X} \psi(x)^2 dx > \left(\int_{x_0}^{X} \varphi(x)\psi(x)dx\right)^2, \qquad (C)$$

for continuous functions of one variable. Although he states it with a strict inequality, Bunyakovsky also describes the case of equality. This inequality is so familiar that it may

---

[1] Université de Reims Champagne-Ardenne, Laboratoire de Mathématiques (CNRS, UMR9008), B.P. 1039, F-51687 Reims Cedex 2. *E-mail* : satyanad.kichenassamy@univ-reims.fr
*Web* : https://www.normalesup.org/~kichenassamy/
[2] For an overview of the history of the Saint Petersburg school of Mathematics, see G.I. Sinkevich and A.I. Nazarov (eds.) (2018). Математический Петербург: История, наука, достопримечательности : Справочник-путеводитель. Санкт-Петербург : Образовательные проекты. English translation forthcoming.
[3] C.A. Andréieff, Викторъ Яковлевичъ Буняковскій. Некрологическій очеркъ, *Communications de la Société mathématique de Kharkow*. 2ème série, (1889), **2**:1-2, 149–161. See also V. E. Prudnikov. *V Y. Bunyakovsky, uchënyj i pedagog* ("V. Y. Bunyakovsky, Scientist and Teacher"; Moscow, 1954) and the comments by Gaïduk (Ю. М. Гайдук, В. Е. Прудников, "В.Я. Буняковский – учёный и педагог" (рецензия), УМН, 1956, том 11, выпуск 2(68), 245–248). Other analyses include A.P. Youschkevitch. *Istoria matematiki v Rossii do 1917 goda* ("History of Mathematics in Russia Before 1917", Moscow, 1968), pp. 296–302 and O.B. Sheynin, "On V. Ya. Buniakovsky's Work in the Theory of Probability", *Arch. Hist. Exact Sciences* **43**(3) (1991) 199-223.
[4] *Mémoires de l'Académie impériale des sciences de St.-Pétersbourg* (VIIᵉ série, tome I, N°9, juillet 1859, 18 pp., read April 29, 1859).

seem almost obvious. For Bunyakovsky, it was a direct consequence of the "well-known inequality"[5]

$$(a_1^2 + \cdots + a_n^2)(b_1^2 + \cdots + b_n^2) > (a_1 b_1 + \cdots + a_n b_n)^2, \tag{1}$$

applied to what we now call "Riemann sums". It is likely that his readers would have recognized a reference to Note II of Cauchy's 1821 lectures, discussed below; in fact, Cauchy's notation is followed in many places of the paper.

In this paper, we try and understand why the Cauchy-Bunyakovsky inequality for integrals was obtained by Bunyakovsky but missed by other authors of the same period who had obtained closely related results[6]. These now familiar results may be viewed anachronistically as consequences of the fact that orthogonal projection on a line, in $\mathbb{R}^n$ endowed with the usual scalar product, decreases the norm. They admit of many generalizations, such as the "Rogers-Hölder inequality". Despite these reinterpretations, modern Russian usage does not entirely conflate all versions of the inequality, distinguishing the Cauchy-Bunyakovsky[7] inequality for integrals, and the abstract Cauchy-Schwarz inequality[8]. This suggests that the difference in terminology reflects a mathematical fact: the same integral inequality may be embedded in several, conceptually different discourses, leading to different proofs and different generalizations. Other traditions are less precise about the identification of theorems. While Bunyakovsky's priority was recognized clearly by Hardy, Littlewood and Pólya[9], they still called his inequality "the inequality of Schwarz" (see the title of section 6.5, p. 135), even though they give two entirely different proofs, suggesting that two different conceptual frameworks are involved. A textbook mention may be found in to 1912, by which time the name "Schwarz' inequality" had apparently become standard[10]. With the emergence of inequalities as an identified research topic, earlier results were revisited, and the name of Cauchy added to Schwarz'. In the process, Bunyakovsky's priority was recognized, even if this fact is not always reflected in the terminology. More important, Schwarz' or Bunyakovsky's motivations were lost, and the results conflated into one.

---

[5] Bunyakovski actually writes $a_1^2 + a_2^2 + a_3^2 + \cdots + a_n^2$, and similarly for the other two sums.
[6] This paper is an outgrowth of a talk delivered at the Saint-Petersburg Seminar on the History of Mathematics on April 15, 2021. We take this opportunity to thank Galina I. Sinkevich for her kind invitation and her remarks.
[7] W. Stekloff (V. Steklov) (« Sur la théorie de fermeture », *Bulletin de l'Académie Impériale des Sciences de St.-Pétersbourg*, VI$^e$ série, 1916, vol. 10(4), 219–226) uses « Bouniakovski's well-known inequality » on p. 225. Five years later, he calls it « l'inégalité de Schwarz-Bouniakowsky » (« Une contribution nouvelle au problème du développement des fonctions arbitraires en séries de polynômes de Tchébychef », *Bulletin de l'Académie des Sciences de Russie*, VI$^e$ série, 1921, vol. 15, 267–280, see p. 275). N. M. Kryloff (Krylov) ("Application of the method of W. Ritz to a system of differential equations", *Bulletin de l'Académie des Sciences*, 11(8) (1917), 521–534) prefers "Bouniakowsky-Schwarz" (p. 530).
[8] Лев Дмитриевич Кудрявцев, *Курс Математического Анализа*, (2006) Т. 3, p. 192 and pp. 188-9.
[9] *Inequalities,* Cambridge, 1934, 2$^{nd}$ ed., 1952, see th. 7, note a, p. 16, and th. 181, pp. 132-4. Hardy and Littlewood had previously called the discrete inequality the "familiar inequality of Cauchy and Schwarz", see p. 38 of G.H. Hardy & J.E. Littlewood, "Some problems of 'Partitio Numerorum' ; I : A new solution of Waring's problem" (*Nachrichten von der Gesellschaft der Wissenschaften zu Göttingen, Mathematisch-Physikalische Klasse*, 1920, Heft 1, presented by E. Landau on January 30, 1920, 33-54). It is applied again on p. 29. The name "inégalité de Schwarz" appears as early as 1887 (H. Poincaré, "La méthode de Neumann et le problème de Dirichlet", *Acta Mathematica* 20 (1897) 59-142 (printed Nov. 11, 1895)). Poincaré calls the discrete inequality "inégalité de Lagrange" (p. 72), from which he derives (like Bunyakovsky) the integral inequality that he calls "inégalité de Schwarz" (p. 73). Like Schwarz, he only writes strict inequalities.
[10] H.B. Heywood, Maurice Fréchet, *L'équation de Fredholm et ses applications à la physique mathématique*, Hermann, Paris, 1912, with a préface and a note by Jacques Hadamard. See §6 of the introduction for "Schwarz' inequality", with the usual proof involving a second-degree polynomial.

We suggest in this paper that Bunyakovsky's inequality should be understood against the backdrop of his mathematical discourse, and that the same goes for each of the inequalities loosely associated with the names "Cauchy-Bunyakovsky-Schwarz". We also show that some essential elements in Bunyakovsky's, or other mathematicians' discourses, such as proofs, are not spelled out, but are implied by their discourses. That makes them *apodictic discourses*, or discourses that contain proof (*apodeixis*) encoded in part by the discursive structure itself: if the motivation is clear, and the formulation precise enough, proofs need not be spelled out[11].

We first examine the inequalities obtained by Cauchy, Liouville, Grassmann, who all failed to identify Bunyakovsky's inequality. This will show that the Cauchy-Bunyakovsky inequality wasn't obvious at the time and that it had conceptually different proofs. We then give an outline of Bunyakovsky's treatment of the inequality, putting his result against the backdrop of the mathematical discourse in which it is embedded. We also briefly analyze later work by Schwarz, Gram, Rogers and Hölder and others, showing how the inequality was reinterpreted and finally incorporated in elementary textbooks with, however, different proofs, and generalizations in various directions.

## 2. What Cauchy, Liouville and Grassmann did, and did not do.

For Cauchy, his inequality is a mere illustration in a discourse[12] on operations with the symbols of inequality $<$ and $>$, applied to the comparison of means. Even though these signs are already found in Harriott's *Artis analyticae praxis* (1631)[13], and although Cauchy is lecturing to advanced students at the École Polytechnique, it is striking that he felt the need to spell out rules of manipulation of inequalities. His inequality[14] (eq. (30) in Th. XVI, p. 373) is for him a consequence of the identity (eq. 31, p. 374)[15]

$$(a_1 b_1 + \cdots + a_n b_n)^2 + \sum_{i<j}(a_i b_j - a_j b_i)^2 = (a_1^2 + \cdots + a_n^2)(b_1^2 + \cdots + b_n^2). \qquad (2)$$

He first obtained the latter in the case when all the $b_j = 1$, that he proved by expanding the squares, and spells out the cases $n = 2$ and $3$, before proving them. He then proves the general result, again by expansion. As special cases (eq. (33), in Scolie II) he recovers: (i) for $n = 2$, an identity he had already used (eq. (8) of chap. I, §8) and (ii) for $n = 3$, an inequality stated to be useful "in the theory of radii of curvature of curves traced on arbitrary surfaces, and in several questions of Mechanics." It had been given by Lagrange in 1773 in connection with Mechanics[16]. In modern terms, these identities read

---

[11] For an introduction to apodictic discourse, see S. Kichenassamy, "Brahmagupta's apodictic discourse", *Gaṇita Bharatī,* **41**(1) (2019) 93-113.

[12] Note II (« Sur les formules qui résultent de l'emploi du signe > ou <, et sur les moyennes entre plusieurs quantités », pp 360-377) of the *Cours d'Analyse de l'École Royale Polytechnique*, (1821) [tome III of *Œuvres Complètes d'Augustin Cauchy*, sér. 2, Académie des Sciences, Paris, 1882-1974].

[13] For the signs of inequality, see Florian Cajori's monograph (*A History of Mathematical Notations*, vol. 1 (Open Court, London, 1928)), §188, p. 199. See also his *A History of Mathematics* (AMS, 2000 [First ed., 1893]), p. 157 : "The signs [of inequality in the wide sense] were first used about a century later [than Harriott] by the Parisian hydrographer Pierre Bouguer" with a reference to P. H. Fuss, *Corresp. Math. Phys.,* I, 1843, p. 304; and to *Encyclopédie des sc. math.*, T. I, vol. I, 1904, p. 23.

[14] Written for square roots.

[15] Cauchy does not use an index notation. The two sets of numbers are $a, a', a'', ...$ and $\alpha, \alpha', \alpha'', ...$, and sums are indicated by their first few terms, as in $a^2 + a'^2 + a''^2 + \cdots$, or $(a\alpha' - a'\alpha)^2 + (a\alpha'' - a''\alpha)^2 + \cdots + (a'\alpha'' - a''\alpha')^2 + \cdots$

[16] J. L. Lagrange (1773), "Nouvelle solution du problème du mouvement de rotation d'un corps de figure quelconque qui n'est animé par aucune force accélératrice," *Nouv. Mém. de L'Acad. Roy. de Berlin,* 85-128 [Œuvres, iii. pp. 577-616]. See also M. Gidea, C.P. Niculescu, "A Brief Account on Lagrange's Algebraic Identity", *Math Intelligencer* **34,** 55–61 (2012).

$$(ac + bd)^2 + (ad - bc)^2 = (a^2 + b^2)(c^2 + d^2) \tag{3}$$

and, in three dimensions,

$$(u|v)^2 + \|u \times v\|^2 = \|u\|^2 \|v\|^2. \tag{4}$$

These in turn admit two natural generalizations in higher dimensions: one is to view all of them as versions of the identity $\cos^2 \theta + \sin^2 \theta = 1$, which leads to (1) in $n$-space and, ultimately, to the abstract Cauchy-Schwarz inequality.[17]

Cauchy states his general inequality again in a "Note sur la détermination approximative des racines d'une équation algébrique ou transcendante[18]" that concludes his lectures on differential calculus form 1829 (see Scolie II, p. 689), and applies it to prove that the modulus of the sum of two complex numbers lies between the sum and the difference or their moduli (in modern notation, $|z| - |z'| \leq |z + z'| \leq |z| + |z'|$).

Identities implying the inequality in low dimensions are attested a little before Cauchy's Note. For n = 4, the inequality is also a consequence of the identity Euler (1778) used in relation to the decomposition of a number in a sum of four squares. Thus, inequality (1) by itself does not have a clear mathematical meaning: is it a consequence of identity (3) or (4)? And what is the mathematical meaning of these identities? For Cauchy, they belong to the same circle of ideas as the inequality between arithmetic and geometric means[19]. Indeed, the inequality $\sqrt{xy} \leq \frac{1}{2}(x+y)$ follows from the identity $(a+c)^2 + (a-c)^2 = 4a^2c^2$, with $x = a^2$ and $y = c^2$, which is a special case of (3) (with $b = d = 1$). And indeed, Cauchy gives this derivation and obtains from it the inequality $\sqrt[n]{a_1 \cdots a_n} \leq \frac{1}{n}(a_1 + \cdots + a_n)$ (Th. XVII, pp. 375-377), that concludes his Note[20].

By contrast, identity (4) does not generalize in a simple way to higher dimensions. For $n = 8$, a similar identity was given by C. F. Degen (1818) in the Mémoires of the St.-

---

[17] The Binet-Cauchy identity is another general result that implies Cauchy's inequality. For the history of this result see pages 80 to 131 of T. Muir, *The theory of determinants in the historical order of its development*, Dover, New York, 1960. Cauchy and Binet presented their results to the French Academy on the very same day (Nov. 30, 1812). They were published in the *Journal de l'École Polytechnique*, **9** (1813) 280-354 (Binet) and **10** (1815) 29-112 (Cauchy).

[18] In *Leçons sur le calcul différentiel*, 1829, Œuvres, sér. 2, t. IV. The note is found on pages 573-609.

[19] He may have been influenced by Maclaurin's papers on inequalities, namely "IV. A letter from Mr. Colin Mac Laurin, Professor of Mathematicks at Edinburgh, and F. R. S. to Martin Folkes, Esq; V. Pr. R. S. concerning Æquations with impossible roots", *Phil. Trans. R. Soc.*, **34** (1727), 104–112, and "IV. A second letter from Mr. Colin McLaurin, Professor of Mathematicks in the University of Edinburgh and F. R. S. to Martin Folkes, Esq; concerning the roots of equations, with the demonstration of other rules in algebra; being the continuation of the letter published in the Philosophical Transactions", N° 394, *Phil. Trans. R. Soc.*, **36** (1730), 59–96 (letter dated Apr. 19, 1729). Indeed, MacLaurin' work seems to have had some influence on the Continent (J. V. Grabiner, "Was Newton's calculus a dead end? The continental influence of Maclaurin's treatise of fluxions," *Amer. Math. Monthly* **104** (1997), 393–410).

[20] He proves it for the case when the number of terms has the form $n = 2^k$, by induction on $k$ – this follows by repeatedly applying the case $k = 1$. Indeed, if the result is known for some value of $n$ terms, we may obtain it for $2n$ terms by writing, $\sqrt[2n]{a_1 \cdots a_{2n}} \leq \frac{1}{2}\left(\sqrt[n]{a_1 \cdots a_n} + \sqrt[n]{a_{n+1} \cdots a_{2n}}\right) \leq \frac{1}{2}\left(\frac{1}{n}(a_1 + \cdots + a_n) + \frac{1}{n}(a_{n+1} + \cdots + a_{2n})\right) = \frac{1}{2n}(a_1 + \cdots + a_{2n})$, as desired. Cauchy then derives the general case from this by padding that is, by considering sequences of the form $a_1, \ldots, a_n, M, \ldots, M$ where $M = \frac{1}{n}(a_1 + \cdots + a_n)$, where $M$ is repeated as many times as necessary so that the length of the sequence is a power of two.

Petersburg Academy[21]. There are identities[22] for other powers of 2. However, they do not, and could not, involve bilinear expressions (Hurwitz, 1898).

Liouville[23] (1836) gives an inequality for the terms of the expansion of a function of one variable in a generalized Fourier series. The inequality he gives on page 265 says, in his notation, that the sum of finitely many expressions of the form $C_n^2$, where $C_n = \int g V_n f dx$, where $g$ is a weight function, is bounded above by $\int g f^2 dx$. Now, by his definition, $C_n = \frac{\int g V_n f dx}{\int g V_n^2 dx}$. Therefore, if the sum in Liouville's expression contains only one term ($n = 1$), we obtain after a short calculation

$$\left(\int g V_1 f dx\right)^2 \leq \int g V_1^2 dx \int g f^2 dx.$$

For $g = 1$, we recover a special case of Bunyakovsky's inequality. As Lützen shows[24], Liouville's discourse does not suggest that he had envisioned the general result that we now call "Bessel's inequality". Not only should he have taken a geometric view and interpreted the left-hand side in terms of the projection of $f$ on the line generated by $V_1$, he should also have considered that any nonzero continuous function can be the first element of what we now call a Hilbert basis. These notions appear much later. Thus, Liouville did not obtain Bunyakovsky's inequality, although he had the mathematical tools to state and prove it, because the line of his mathematical discourse did not point to it.

Grassmann[25] (1862), in his *Ausdehnungslehre*, takes a fully geometric view, but does not apply it to functions. In addition, he does state the abstract form of the "Cauchy-Schwarz" inequality but does not prove it[26]! He defines the length of the "shadow" (orthogonal projection) of one line on another and calls the ratio between the two the cosine of an angle, without checking that the "cosine" thus defined does not exceed 1 in absolute value. In his notation, he writes[27] $\cos \angle AB = \frac{[A|B]}{\alpha\beta}$ with $\angle AB$ between 0 and $\pi$. The proof would have been easy to supply in his formalism: Let us spell out the argument assuming $B \neq 0$. In n°166 (p. 104) Grassmann defines the "shadow" $A'$ of line $A$ by orthogonal projection on the direction of $B$. From this, it follows that $A'$ and $A - A'$ have vanishing scalar product. Hence, the sum of the squares of the lengths of $A'$ and $A - A'$ is equal to the squared length of $A$. The length of $A'$ therefore cannot exceed that of $A$. In modern notation, using parentheses where Grassmann uses square brackets, $A' = (B|A)B/(B|B)$. Since $(A'|B) = (A|B)$ by construction,

---

$(A - A'|B) = 0$. Since $A'$ is a multiple of $B$, $(A - A'|A') = 0$. Therefore, $(A|A) = (A - A' + A'|A - A' + A') = (A - A'|A - A') + (A'|A')$. Hence $(A'|A') \leq (A|A)$, with equality precisely when $A$ is a multiple of $B$. Grassmann may well have had this proof in mind, but does not seem to have separated his geometric intuition from the abstract formalism, so that no proof appeared necessary to him.

Thus, Grassman's discourse is apparently formalized, but relies on a geometric intuition that causes him to fail to identify a crucial point. The emphasis on geometric intuition probably explains why applications to sets of data, or to functions, were not considered.

Let us now turn to Bunyakovsky's line of thought.

### 3. Bunyakovsky's paper : a systematic investigation of means.

Bunyakovsky opens his 1859 paper by proposing that the notion of arithmetic mean is the concept that, when applied to "functions of one or several variables that change by imperceptible steps leads to Integral Calculus, in the most natural, the most elegant and the most satisfactory way from the point of view of clarity" (all translations are ours). In retrospect, his point of view is very close to that of Hardy, Littlewood and Pólya[28], who will similarly explore many years later the various types of means, and their relations. By contrast, modern courses would rather appeal to notions of length or area, that are less general, since they are limited to special types of magnitudes.

Bunyakovsky continues: "In a large number of applications of transcendental Analysis, this point of view greatly facilitates the conception of the relations that exist between the data in the question [at hand], as in many examples that could be mentioned, in the Theory of Probability among others." A footnote gives the only reference of the paper: "See my Traité du Calcul des Probabilités (Основанія Математической Теорія Вѣроятностей [Foundations of the Mathematical Theory of Probability], 1846 г.)".[29] Continuous means appear on p. 446 sq. of this treatise, the discrete ones on p. 153, and the convergence rules of Duhamel and Raabe, mentioned in §4 of the paper appear on p. 397 of his book. Bunyakovsky then introduces a notation for means that is also found in his Treatise (p. 446 sq.).[30]

$$M_{x_0}^X f(x) = \frac{\int_{x_0}^X f(x) dx}{X - x_0}.$$

In many cases, it is implicitly assumed that $f$ is positive. The context makes the appropriate conditions obvious.

The paper proceeds in five paragraphs. §1 derives five inequalities (A–E), and §2 gives the mean-value theorem[31] probably taken from Cauchy, with numerical examples. Indeed, he was thoroughly familiar with Cauchy's work, whose 1823 lectures he had freely translated into

---

[28] *Inequalities*, op. cit.

[29] We have only seen it in Russian; we have not found evidence of a translation. Bunyakovsky also refers to convergence criteria of Duhamel and Raabe in §4, but does not give a precise reference.

[30] In the following, $X$ and $x_0$ should respectively stand above and below the $M$ sign. We weren't able to achieve this effect for typographical reasons.

[31] The history of the mean-value theorem has been considered in many recent works; see e.g. Sinkevich (« Rolle's theorem and Bolzano-Cauchy theorem : A view from the end of the 17th century until K. Weierstrass' epoch», *Gaṇita Bhāratī*, **38**:1 (2016) 31-53), J. Barrow-Green ("From cascades to calculus: Rolle's Theorem", in *The Oxford Handbook of the History of Mathematics* (eds. Eleanor Robson, Jacqueline Stedall), Oxford University Press, 2009, 737- 754) and C. Smoryński (*MVT : A Most Valuable Theorem*, Springer, 2017).

Russian with his comments.[32] Inequality (C) is "Bunyakovsky's inequality". The results are then applied in the remaining sections. §3 derives inequalities on remarkable numbers; §4 is devoted to the convergence of special series; and §5 gives discrete analogues (A'–E') of results (A–E). Throughout the paper, he uses < or > where we would often write ≤ and ≥; he is probably following Cauchy in this. Still, he seems quite aware of the correct conditions for equality. Let us give a sample of his results in §§1–3, focusing on those that involve (C).

In §1, Bunyakovsky first extends geometric and harmonic means to the continuous case: the geometric mean of a function is the exponential of the arithmetic mean of its logarithm. Its harmonic mean is the reciprocal of the arithmetic mean of $\frac{1}{f}$.

(A): the arithmetic mean exceeds the geometric mean.
(B): the geometric mean exceeds the harmonic mean. Applying this to the logarithm of f yields a two-sided bound for the average of $f$. This is the main tool for §2.
(C): Bunyakovsky's inequality:
$$\int_{x_0}^{X} \varphi(x)^2 dx \int_{x_0}^{X} \psi(x)^2 dx > \left( \int_{x_0}^{X} \varphi(x)\psi(x) dx \right)^2.$$

(D): Special case of (C) where $\varphi = \frac{1}{\psi} = \sqrt{f}$ :
$$\int_{x_0}^{X} f(x) dx \int_{x_0}^{X} \frac{1}{f(x)} dx > (X - x_0)^2 ;$$

(E): Special case of (C) where $\psi = 1$:
$$\left( \int_{x_0}^{X} \varphi(x) dx \right)^2 < (X - x_0) \int_{x_0}^{X} \varphi(x)^2 dx.$$

(C) is for him a consequence of Cauchy's inequality (1); he adds that equality holds for $\psi(x) = \lambda \varphi(x)$, where $\lambda$ is a constant; he notes a few lines below, that there is also equality for $X = x_0$. Even though Bunyakovsky does not supply a proof, the discussion of equality is a consequence of the integral form of the identity from which Cauchy had derived his own inequality. This identity is given by Hardy, Littlewood and Pólya[33], and Courant and Hilbert[34]: with obvious notation,

$$(f, g)^2 = (f, f)(g, g) - \frac{1}{2} \iint (f(x)g(y) - f(y)g(x))^2 dx dy.$$

In §2, Bunyakovsky derives the result $F(X) = F(x_0) + (X - x_0)F'(x_0 + \lambda(X - x_0))$, for functions $F$ of class $C^1$, using first his observation that (A) and (B) applied to $\log f$ with $f = F'$ leads to a two-sided bound on the average of $f$. He then observes that the two sides involve the exponentials of $F'$ and $-F'$. He also points out in a footnote that the result could be obtained directly from the mean-value theorem for integrals, applied to $F$. In some examples, the value of $\lambda$ may be obtained explicitly.

§3 lists a number of applications of (A-D). They include six different applications of (D) [a special case of (C)], giving inequalities such as $\pi^2 + 2\pi > 16$, or $9 \ln 3 > \pi^2$, where we

---

[32] We thank Galina I. Sinkevich for this remark. A Russian translation of Cauchy's 1821 lectures appeared in 1864.
[33] *Inequalities,* (op. cit.) p. 133.
[34] *Methods of Mathematical Physics*, Wiley-Interscience, 1989, p. 49 [Interscience, 1953; Springer, 1937].

corrected a misprint in paper [35] – the former means that $16,15... > 16$, and the latter, $9,8874... > 9,8696...$. To prove the first of these inequalities, take $f(x) = \sqrt{1-x^2}$, $x_0 = 0$ and $X = 1/\sqrt{2}$; for the other, take $f = \sin x \cos x$, $x_0 = \pi/6$, $X = \pi/4$.

## 4. Gram (1882) and Schwarz (1885): quadratic forms.

Schwarz' 1885 paper[36] seems to be the source of the most commonly taught proof, relying (in modern notation) on the fact that, if a norm is associated with a real scalar product, then $p(t) := \|x + ty\|^2$ is a second-degree polynomial in $t$ and therefore, has nonpositive discriminant. However, Schwarz' argument is slightly different and reveals the motivation behind his proof. Schwarz (§15, p. 344) derives (C) in two variables – Bunyakovsky's $\psi$ is now called $\chi$ – from the fact that $\iint (\alpha\varphi + \beta\chi)^2 \, dx \, dy$ is a positive-definite binary quadratic form in $\alpha$ and $\beta$ and therefore, that its discriminant is positive, assuming, as he does, that the quotient of the two functions is not identically constant, and that the integrals of $\varphi^2, \varphi\chi$ and $\chi^2$ are all absolutely convergent. The minimum of the integral is of course positive unless the two functions are proportional to one another, hence his assumption; like Cauchy and, later, Poincaré, Schwarz uses only the sign for strict inequality. The introduction of a minimum problem may have been motivated by Weierstrass' lectures on the calculus of variations, of which there are at least ten manuscript notes by students between 1864 and 1882. To justify his result, Schwarz appeals to the "known result that the discriminant of a positive binary quadratic form always has a positive value, different from zero". The proof and motivation are only indicated by the choice of words: he does not refer to the discriminant of a quadratic *polynomial*, as in a popular proof, but to the discriminant of a *binary form*. This is a transparent allusion to the invariant theory, which he very likely learnt from Clebsch's book[37]. Clebsch had made quite a sensation in Leipzig, where he taught and passed away shortly before Schwarz was appointed in Leipzig. Clebsch had already envisioned geometric interpretations of invariant theory (*op. cit.* §§17-25).

Had Schwarz interpreted the discriminant as a determinant, he would have been led to the consideration of the more general quadratic form

$$\|\alpha_1\varphi_1 + \cdots + \alpha_n\varphi_n\|^2 = \iint (\alpha_1\varphi_1 + \cdots + \alpha_n\varphi_n)^2 \, dx \, dy, \qquad (Q)$$

and therefore, to the positivity of the determinant of Gram matrices[38]. The latter had appeared in a recently published paper (Gram[39], 1883, p. 43). Here again, the discursive structure into which Bunyakovsky's inequality is embedded – the choice of a single word, "discriminant" over "determinant"! – determines possible generalizations. Such a generalization is pointed out by Richardson and Hurwitz (1909)[40], who refer to E. Schmidt's (1907, 1908) for the geometric interpretation, and mention the relation to the method of least squares, but without reference to Gram (1883). For the case of the trigonometric system, the relation to least-squares seems to go back to Bessel (1828)[41] and to Gauss' theory of errors[42]. By contrast, the usual argument, by referring only to the one-variable polynomial $p(t) = \iint (\varphi + t\chi)^2 \, dx \, dy$, does not suggest the introduction of Gram determinants.

---

[35] One reads $9\sqrt{3}$ instead, which is much too large.
[36] « Über ein die Flächen kleinsten Flächeninhalts betreffendes Problem der Variationsrechnung », *Acta Societatis Scientiarum Fennicæ*, T. 15 (1885) 315-362. Presented to Karl Weierstrass on Oct. 31, 1885.
[37] A. Clebsch, *Theorie der binären algebraischen Formen,* Teubner, Leipzig, 1872, see §33.
[38] The Gram matrix of the set $(\varphi_1, ..., \varphi_n)$ is the matrix of the quadratic form (Q) in the variables $\alpha_1, ..., \alpha_n$. It is the square matrix made up with the scalar products $(\varphi_i | \varphi_j)$ $(1 \leq i, j \leq n)$. If the set is independent, the form (Q) is positive-definite and therefore, the determinants of the matrices $\left((\varphi_i|\varphi_j)\right)_{1 \leq i,j \leq k}$ for $k = 1, 2, ..., n$ are all positive. For $k = 2$, this recovers the Cauchy-Schwarz inequality $(\varphi_1|\varphi_2)^2 \leq (\varphi_1|\varphi_1)(\varphi_2|\varphi_2) = \|\varphi_1\|^2 \|\varphi_2\|^2$.
[39] J. P. Gram, "Ueber die Entwickelung reeller Functionen in Reihen mittelst der Methode der kleinsten Quadrate", *Journal für die reine und angewandte Mathematik*, 94(1) (1883), 41 – 73.
[40] R. G. D. Richardson, W. A. Hurwitz, "Note on determinants whose terms are certain integrals", *Bull. Amer. Math. Soc.* 16(1): 14-19 (October 1909).

### 5. Generalizations

The geometric interpretation of the integrals of the previous paragraph seems to have been appreciated by M. Fréchet (1907)[43] and E. Schmidt. Once this step had been taken, the interpretation of the inequality in terms of the contraction property of orthogonal projection becomes natural. However, they are so diverse that it does not seem possible to subsume all of them under a single overarching theory. The question "What *is* the mathematical meaning of the Cauchy-Bunyakovsky inequality" does not seem to be meaningful, insofar as the inequality has several interpretations[44]. Let us review some striking ones.

The very useful fact that the inequality remains valid for Hermitian forms $q(x, y)$, nonnegative but not necessarily positive definite, seems to have been folklore in the 1930s and possibly earlier. It is proved in a special setting by von Neumann[45], and by Riesz and Nagy[46], who do not seem to consider it as new. It is this general form that Bourbaki calls the "Cauchy-Schwarz inequality" (*EVT*, V.3, Prop. 2). Their proof is in two steps. If $q(y, y) \neq 0$, derive the result from $q(y, y) q(x + \xi y, x + \xi y) \geq 0$, by inserting $\xi = -\overline{q(x, y)}/q(y, y)$. This choice is natural from the point of view of projections, since $-\xi y$ is the the projection of $x$ on the line generated by $y$. If $q(y, y) = 0$, they note that $q(x + \xi y, x + \xi y) = q(x, x) + \xi q(x, y) + \overline{\xi q(x, y)}$ must be nonnegative. But a linear form must change sign unless it vanishes identically. Specifically, Bourbaki sets $\xi = -\overline{q(x, y)}$ to obtain $q(x, y) = 0$. Thus, in all cases, $|q(x, y)|^2 \leq q(x, x) q(y, y)$. This argument is hybrid: the first part is motivated by projection, and the second by the theory of forms.

The Hölder-Rogers inequality[47] is a completely different generalization of Cauchy's inequality that relies on convexity. Von Neumann's proof of the abstract inequality in Hilbert space[48] is actually reminiscent of the usual proof of Hölder's inequality. He also proves the analogue for Hermitian forms. This line of thought has been further generalized. Zarantonello[49] gives a necessary and sufficient condition for a set-valued operator $T$ from a real, locally convex topological vector space $X$ with values in (subsets of) its dual $X^*$ to satisfy an inequality of the form $\langle Ty, x \rangle \leq \langle Tx, x \rangle^{\frac{1}{2}} \langle Ty, y \rangle^{\frac{1}{2}}$, or more generally, if $T$ is set-valued, $\langle y^*, x \rangle \leq \langle x^*, x \rangle^{\frac{1}{2}} \langle y^*, y \rangle^{\frac{1}{2}}$ for all $x$ and $y$ in the domain $D(T) \coloneqq \{x \in X \mid Tx \neq \emptyset\}$ of the map $T$, and $x^* \in Tx$, $y^* \in Ty$. The condition is that $Tx \subset \partial k(x)$, for all $x \in D(T)$, where $k$ is a lower semicontinuous, positively homogeneous convex[50] function. In this case, $k$ has the same domain as $T$, and is necessarily given by $k(x) = \frac{1}{2}\langle x^*, x \rangle$, for all $x \in D(T)$ and all $x^* \in Tx$. While this result is noteworthy as it yields a constraint on the functional set-up for the inequality to hold, it is not clear that it includes all possible generalizations: to the best of our knowledge, a proof that all the known generalizations of the Cauchy-Bunyakovsky inequality given in the literature may indeed be cast in this form has not been carried out.

---

[41] "Ueber die Bestimmung des Gesetzes einer periodischen Erscheinung", *Astronomische Nachrichten,* No. 136, 333-348; see equation (4) on p. 340.
[42] Bessel takes his notation from Gauss (*loc. cit.*, previous note).
[43] They seem to have suggested a geometric approach in function spaces; according to Fréchet (*Comptes Rendus de l'Académie des Sciences de Paris*, séance du 20 mars 1905, pp. 772-774), Weierstrass' lectures on the Calculus of Variations suggested to Fréchet the introduction of the notions of neighborhood and distance.
[44] Similar considerations could be developed about all inequalities. In addition to the volume by Hardy, Littlewood and Pólya already quoted, surveys of inequalities include E. F. Beckenbach and R. Bellman (1961). *Inequalities*, Springer; D. S. Mitrinović, J. E. Pečarić, and A. M. Fink (1993), *Classical and New Inequalities in Analysis*, Kluwer, Dordrecht. See also D. S. Mitrinović (1970). *Analytic Inequalities*, Springer and P. S. Bullen, D. S. Mitrinović, and P. M. Vasić (1988) *Means and Their Inequalities*, Kluwer.
[45] *Les fondements mathématiques de la Mécanique quantique,* (tr. A. Proca), Félix Alcan, Paris, 1946, Th. 19 in §5 of chapter II, p.74 (German : *Mathematische Grundlagen der Quantenmechanik*, Springer, 1932, p. 53).
[46] Both work with the expression $(Af, g)$, where $(f, g)$ is a Hermitian scalar product in a Hilbert space, and $A$ is a bounded self-adjoint operator (Frédéric Riesz and Béla Sz.-Nagy, *Leçons d'Analyse Fonctionnelle,* 2ᵉ éd., Akadémiai Kiadó, Budapest, 1953, 104, p. 260). They call it « inégalité de Cauchy-Schwarz généralisée ». They also prove Bunyakovsky's inequality by the same method as von Neumann (§21, p. 41), referring to Hardy, Littlewood and Pólya for the attribution. The usual proof is given for abstract spaces in §83, p. 196-197.
[47] The inequality now known as Hölder's inequality was first given by Rogers. The one known as Jensen was first given by Hölder, as Jensen points out in his paper. See L. J. Rogers, "An extension of a certain theorem in inequalities", *Messenger of Math.*, **17** (1888) 145-150 ; O. Hölder, "Ueber einen Mittelwerthssatz", *Nachr. von der Königl. Gesellschaft der Wiss. und der Georg-August-Universität zu Göttingen*, (1889) 38-47; and J.LW.V.

T.S. Motzkin proposes a different interpretation of the inequality. He argues[51] that "[m]any matrix inequalities concern signs of minors; the two well-known ones that do not, the C-B-S inequality and Hadamard's inequality, are immediate consequences of one that does" As far as Cauchy's inequality is concerned, his point is that Cauchy's inequality expresses that the quadratic form $(\sum_k a_k^2)(\sum_k x_k^2) - (\sum_k a_k x_k)^2$ (in the variables $x_k$) is positive semi-definite, and it is well-known that this positivity can be tested using the signs of the principal minors of the matrix of the quadratic form.

There are also non-commutative operator-theoretic generalizations of Hölder's inequality[52]. Countless other generalizations have been proposed[53].

---

Jensen, "Sur les fonctions convexes et les inégalités entre les valeurs moyennes", *Acta Mathematica,* **30**, (1906) 175-193. Jensen notes that "his basic formula (5) is not entirely new", referring to Hölder's paper (p. 192), but stresses that his smoothness assumptions are weaker.

[48] *Fondements mathématiques de la Mécanique quantique,* op. cit., p. 28-29 (original, p. 21-22).

[49] Eduardo H. Zarantonello, "The meaning of the Cauchy-Schwarz Buniakovsky inequality", Proc. AMS, **59**:1 (1976) 133-137.

[50] Recall that a convex $k$ function takes its values in $(-\infty, +\infty]$, upper bound included, and that its domain is, by definition, the set of $x \in X$ for which $k(x) \neq +\infty$.

[51] T.K. Motzkin, "Signs of minors", in *Inequalities* (Oved Shisha, ed., Academic Press, London, 1967, pp. 225-240, see p. 225). This refers to his Th. 16, pages 237-8. For generalizations in a slightly different direction, see M. Marcus, "The Cauchy-Schwarz inequality in the exterior algebra", *Quarterly J. of Math.* **17** (1966), 61-63, and Antonino M. Sommariva, "The generating identity of Cauchy-Schwarz-Bunyakovsky inequality", *Elemente der Mathematik,* **63** (2008) 1-5. The Cauchy-Binet identity, already mentioned, belongs to the same circle of ideas.

[52] See e.g. Masatoshi Fujii, Saichi Izumino, Ritsuo Nakamoto and Yuki Seo ("Operator inequalities related to the Cauchy-Schwarz and Hölder-McCarthy inequalities", Nihonkai Math. J., **8** (1997), 117-122) and its references.

[53] Without any claim at completeness, see for instance, in addition to the survey already mentioned, H. and N. Sedrakyan (2018). "The Cauchy-Bunyakovsky-Schwarz Inequality" in: *Algebraic Inequalities. Problem Books in Mathematics*; S. M. Sitnik (2010), "Generalized Young and Cauchy-Bunyakovsky inequalities: A survey" in *Advances in Modern Analysis* (Yu. F. Korobeinik and A.G.Kusraev eds.), South Mathematical Institute of the Vladikavkaz Scientific Center of the Russian Acad. Sci, and the Government of Republic of North Ossetia–Alania, Vladikavkaz, 2009, pp. 221–266. https://arxiv.org/abs/1012.3864, see also P. Agarwal, A.A.Korenovskii and S. M. Sitnik, in: *Trends in Mathematics. Advances in Mathematical Inequalities and Applications* (eds. P. Agarwal, S.S. Dragomir, M. Jleli, B. Samet), Birkhäuser/Springer), 2018, Chapter 18, 333-349. The last two papers stress that the variables in Hölder's inequality do not play symmetric roles, and deduces an improvement from this observation. For unusual inequalities, see Alawiah Ibrahim and Silvestru Sever Dragomir, « A survey of Cauchy-Bunyakovky-Schwarz inequalities for power series », in G.V. Milovanović and M.Th. Rassias (eds.), *Analytic Number Theory, Approximation Theory, and Special Functions. In honor of Hari M. Srivastava*, Springer, 2014, p. 244-295, and its references. For more applications of the Cauchy-Bunyakovsky inequality and some related ones, see also J. Michael Steele, *The Cauchy-Schwarz Master Class,* Cambridge U. Press, 2004.

Generalization may also be hindered by preconceptions. Hardy et al.[54] had considered that Hölder's inequality could not admit of an algebraic proof, unlike Cauchy's inequality: "Cauchy's inequality […] is a proposition of finite algebra […]. It is a recognised principle that the proof of such a theorem should involve only the methods of the theory to which it belongs. […the status of] theorems, such as Hölder's inequality […] depends upon the value of a parameter $k$. If $k$ is irrational […] it is obvious that there can be no strictly algebraical proof." Nonetheless the (Rogers-)Hölder inequality may be obtained by means of *algebraic identities*, and this leads to new, optimal improvements of this inequality[55]. In other words, if $k$ is treated as a variable, then the problem does have an algebraic structure. The main point is the existence of functional equations satisfied by the $R(\alpha, a, b) := \alpha a + \beta b - a^\alpha b^\beta$ where $\beta = 1 - \alpha$. It admits an exact, manifestly nonnegative expression for rational $\alpha$, but also exact addition and multiplication theorems relating its values for different exponents $\alpha$. Earlier refinements of Hölder's inequality are consequences of the computation of $R$ for $\alpha = \frac{1}{2}$. Using these functional equations, one can estimate the difference between the value of the remainder function $R$ for any $\alpha$ and its value for any rational exponent $r$. This difference may be made arbitrarily small, so that our results are optimal.

Thus, the generalizations of Bunyakovsky's inequality have not been subsumed under a single "master theorem" in any natural way. The reason seems to be that they may be derived from *several* non-equivalent identities. Therefore, replacing identities by inequalities introduces an element of indetermination. Each generalization is, in fact, an interpretation that, in turn, suggests new results, to be subjected to further interpretations.

### 6. Conclusions

Bunyakovsky found his inequality because he was interested in means that occur in Statistics and Probability, from a mathematical perspective inherited from Cauchy. This seems to have been a unique combination at the time. By contrast, Schwarz was motivated by Weierstrass' view of the Calculus of Variations and by the Theory of Invariants as expounded by Clebsch. The inequality was first called after Schwarz, whose paper suggested the proof most often given today; however, his actual proof was slightly different in spirit. If his paper had been properly analyzed, it would have immediately led to the positivity of Gram determinants. Only gradually was the complexity of the historical development realized, when inequalities became the focus of interest in their own right.

Later views are reinterpretations of these earlier results, leading to theories of orthogonal systems, from the point of view of analysis (Liouville), data analysis and statistics (Bessel, Gram) or geometry and what we now call exterior algebra (Grassmann), among others. Interpretations of Bunyakovsky's inequality include the following.

- It is the consequence of an identity[56].
- It follows from the positivity of a quadratic form.
- It follows from the convexity of a quadratic form.
- It expresses that orthogonal projection shortens lengths.
- It expresses the positivity of the squared error in the method of least squares.

---

[54] *Inequalities, op. cit.*, p. 7.
[55] S. Kichenassamy, "Improving Hölder's inequality", *Houston Journal of Mathematics*, **36** (1) (2010), 303-312, https://hal.archives-ouvertes.fr/hal-00826949v1
[56] There are, as we have seen, *several* equalities that lead to the same inequality.

The various proofs of the modern inequality suggest that mathematical tradition has kept a trace of the consciousness that we deal with several related results, but with quite different mathematical meanings. Each interpretation makes some generalizations natural and makes others difficult to guess. What we now call the « Cauchy-Bunyakovsky-Schwarz » inequality is not a single result, but the conflation of several results with different mathematical contexts: projection, convexity, least squares approximation, etc. The various interpretations of Bunyakovsky's inequality are therefore not logical developments of a single idea. Unlike nested Matryoshka dolls, each of which is a copy of the next one, only smaller, the sequence of these inequalities rather resembles Fabergé eggs (Яйца Фаберже), each containing a "surprise[57]" entirely different from the outer shape: the outer (newer) layer does not indicate what the interior (the past) contains. That is why historical investigation provides a deeper understanding of modern mathematics that cannot be inferred from the modern stage alone.

From a wider perspective, this paper suggests that:

- Theorems are best understood against the backdrop of the specific discourse in which they arise.
- Apodictic discourse, that makes motivation and derivations transparent, is conducive to progress by inserting the result in a continuity and thus, suggesting new results.
- In teaching, it seems advisable to label theorems by their meaning, to attach names if necessary to proofs or discourses only, and to give several proofs of the "same" theorem.

---

[57] Fabergé eggs are gifts to the Russian imperial family in Saint-Petersburg by the jeweler Fabergé. From the outside, they appear to be eggs studded with precious stones; when opened, they are found to contain a new object or mechanism called the *surprise*, that would be impossible to guess from the outside. Examples include the recently recovered Third Imperial egg from 1887, the Diamond Treillis from 1892, or the Lilies of the Valley from 1898, that illustrate different types of "surprises". Thus, the 1892 egg contained an elephant. By contrast, the Matryoshka dolls or "Russian dolls", usually shaped in the form of a woman in traditional attire, contain within them a smaller copy of the same doll that, in turn, contains another, smaller copy of itself, and so forth. Similarly there are two points of view of history: in one view, the past is an earlier, cruder stage of a single line of evolution leading to the present, so that its knowledge shall not alter our present views; in the other, the study of the past leads to "surprises" that may help us see the present in a different light, and map out future work.